\documentclass[12pt]{amsart}
\usepackage[centertags]{amsmath}
\usepackage{amsfonts,amssymb,mathrsfs}
\usepackage{colonequals}
\usepackage{graphicx}
\usepackage{placeins}
\usepackage{enumerate}
\usepackage[latin1]{inputenc}
\usepackage{float}

  \newtheorem{theorem}{Theorem}
  \newtheorem{corollary}{Corollary}
  \newtheorem{proposition}{Proposition}
  \newtheorem{lemma}{Lemma}%
  \theoremstyle{remark}

  \newcommand{\NN}{\mathbb N}

 \newcommand{\Mod}[1]{\ (\mathrm{mod}\ #1)}

\begin{document}

\title{Shifted-prime divisors}
\date{\today}

\author{Kai (Steve) Fan}
\address{Max-Planck-Institut f\"{u}r Mathematik, 53111 Bonn, Germany}
\email{steve.fan@mpim-bonn.mpg.de}
\author{Carl Pomerance}
\address{Mathematics Department, Dartmouth College, Hanover, NH 03755, USA}
\email{carlp@math.dartmouth.edu}

\dedicatory{
For Helmut Maier on his seventieth birthday}

\begin{abstract}
Let $\omega^*(n)$ denote the number of divisors of $n$ that are shifted primes,
that is, the number of divisors of $n$ of the form
$p-1$, with $p$ prime. Studied by Prachar in an influential paper from 70 years ago,
the higher moments of $\omega^*(n)$ are still somewhat a mystery.  This paper
addresses these higher moments and considers other related problems.  

\end{abstract}

\maketitle

\section{Introduction}
\label{S:intro}

Let $\omega(n)$ denote the number of different primes that divide $n$.  This function
has been well-studied, and in particular we know that
\begin{equation}
\label{eq:omega}
\begin{aligned}
\frac1x\sum_{n\le x}\omega(n) &=\log\log x+O(1),\\
\frac1x\sum_{n\le x}\omega(n)^2&=(\log\log x)^2+O(\log\log x),
\end{aligned}
\end{equation}
after results of Hardy--Ramanujan and Tur\'an.  Further,
$\omega(n)$ obeys a normal distribution as given by the
Erd\H os--Kac theorem.  For extreme values, we know that
\begin{equation}
\label{eq:omegamax}
\omega(n)\le (1+o(1))\log n/\log\log n,~~n\to\infty,
\end{equation}
a best-possible result of Ramanujan.

Consider the analogous function $\omega^*(n)$ which counts the number of
shifted prime divisors of $n$, that is, the number of divisors of $n$ of the form
$p-1$, with $p$ prime.  One might guess that assertions like \eqref{eq:omega}
and \eqref{eq:omegamax} hold as well for $\omega^*$.  And in fact, it is easy
to prove that
\begin{equation}
\label{eq:omega*1}
\frac1x\sum_{n\le x}\omega^*(n)=\log\log x+O(1).
\end{equation}
However, the analogy stops here.  As it turns out, the function $\omega^*$ is
considerably wilder than $\omega$.  In some sense, $\omega^*$ is closer to
the total number $\tau(n)$ of divisors of $n$.  For example, after work of Prachar \cite{Pr} we
have $\omega^*(n)\ge n^{c/(\log\log n)^2}$ for some positive constant $c$
and infinitely many $n$.  This was improved in \cite[Proposition 10]{APR} to
\begin{equation}
\label{eq:max}
\omega^*(n)\ge n^{c/\log\log n}
\end{equation}
for a positive constant $c$ and infinitely many $n$, a result which is clearly
best possible, but for the choice of $c$, due to the upper bound
\[
\omega^*(n)\le\tau(n)\le n^{(\log 2+o(1))/\log\log n},~~n\to\infty,
\]
a result due to Wigert.  (Also see  \cite[Section 3]{AGP}.)

This paper deals with the moments
\[
M_k(x)\colonequals\frac1x\sum_{n\le x}\omega^*(n)^k,
\]
for $k=2$ and 3.
Prachar \cite{Pr} showed that
\begin{equation}
\label{eq:S1}
M_2(x)\ll(\log x)^2.
\end{equation}
In a letter to the same journal, Erd\H os \cite{EP} proved that 
\begin{equation}
\label{eq:T1}
S_2(x)\colonequals\frac1x\sum_{[p-1,q-1]\le x}1\ll (\log\log x)^3,
\end{equation}
and indicated how the exponent 3 could be replaced by 1, and possibly even by 0.
Here, $p,q$ run over prime numbers and $[a,b]$ denotes the
least common multiple of $a$ and $b$. 
The connection of
these results on $S_2(x)$ to Prachar's theorem is as follows.  We have
\begin{equation}
\label{eq:S}
M_2(x)=\frac1x\sum_{[p-1,q-1]\le x}\left\lfloor\frac x{[p-1,q-1]}\right\rfloor,
\end{equation}
so that \eqref{eq:T1} and a partial summation argument imply that
\begin{equation}
\label{eq:S2}
M_2(x)\ll \log x(\log\log x)^3,
\end{equation}
with the same remark pertaining to the exponent 3.

In a recent paper Murty and Murty \cite{MM}, completed the proof that
\begin{equation}
\label{eq:S3}
M_2(x)\ll\log x,
\end{equation}
and they showed the lower bound
\begin{equation}
\label{SL1}
M_2(x)\gg(\log\log x)^3,
\end{equation}
which improves on the trivial bound 
\[M_2(x)\ge \left(\frac{1}{x}\sum_{n\le x}1\right)^{-1}\left(\frac{1}{x}\sum_{n\le x}\omega^*(n)\right)^2\gg(\log\log x)^2\]
implied by (\ref{eq:omega*1}) and the Cauchy--Schwarz inequality. Further, they made the conjecture that there is a positive constant $C$ such that
\begin{equation}
\label{eq:SC}
M_2(x)\sim C\log x,~~x\to\infty.
\end{equation}

The topic was picked up again by Ding \cite{D1} who, using the claim 
\begin{equation}
	\label{eq:mm0}
\sum_{p,q\le x}\frac1{[p-1,q-1]}=\sum_{[p-1,q-1]\le x}\frac1{[p-1,q-1]}+O(1)
\end{equation}
from \cite[Equation (4.8)]{MM} (also see \eqref{eq:mm} in Section \ref{S:tail}), showed
that 
\begin{equation}
\label{eq:D}
M_2(x)\gg\log x.
\end{equation}
Further in \cite{D2}, Ding, Guo, and Zhang gave a heuristic argument for the Murty--Murty conjecture \eqref{eq:SC} based on the Elliott--Halberstam conjecture,
with $C=2\zeta(2)\zeta(3)/\zeta(6)\approx 3.88719$.

However, as it turns out, there is an error in the proof of \eqref{eq:mm0}.  In particular,
the error term $O(x)$ there, which results from removing the floor symbol in \eqref{eq:S},
is only valid for those pairs $p,q$ with $[p-1,q-1]\le x$ and not for all pairs $p,q\le x$.
We show below in Section~\ref{S:ding} how a modified version of  Ding's argument \cite{D1}  can save the proof of
\eqref{eq:D}.  Further, we show that not only is the proof of \eqref{eq:mm0} in
error, but the assertion is false, see Section \ref{S:tail}.  This unfortunately seems to invalidate the heuristic in
\cite{D2}.  We certainly agree that the Murty--Murty conjecture \eqref{eq:SC} holds, but
we think the correct constant is closer to 3.1.  We give the results of some calculations
that support this.

We conjecture that $M_k(x)\asymp(\log x)^{2^k-k-1}$ and prove this for the third moment
$M_3(x)$.
The proof is considerably more involved than the second moment, but hopefully we have
presented it in a manner that leaves open the possibility of getting analogous results for
higher moments.

We also consider the level sets $\{n:\omega^*(n)=j\}$, showing that for each fixed 
positive integer $j$,
the natural density exists and is positive, with the sum of these densities being 1.

It may be worth pointing out that our methods used to treat the moments of
 $\omega^*(n)$ can be used
to deal with the natural generalization where $p-1$ is replaced with $p+a$ for a
fixed integer $a\ne0$.

Throughout we let $p,q,r,s,\ell$ run over prime numbers.  We let $(m,n)$ denote the
greatest common divisor of $m,n$, and $[m,n]$ their least common multiple.  We also use
the standard order notations $\ll,\asymp,\gg$ from analytic number theory.

\section{The constant $C$ in (\ref{eq:SC})}
\label{S:const}
In Section \ref{S:tail} we shall prove
Theorem \ref{th:lowerbound} which not only shows that the correction that we make to Ding's proof of the lower bound for $M_2(x)$ is necessary (see Section \ref{S:ding}), but it also suggests that the constant $C=2\zeta(2)\zeta(3)/\zeta(6)\approx 3.8872$ for the Murty--Murty conjecture shown by the heuristic argument given in \cite{D2} is probably incorrect. So, what is the correct value of $C$? We leave this
as an unsolved problem, but perhaps it is helpful to look at some actual numbers.
We have numerical calculations of the values of $M_2(x)=\frac1x\sum_{n\le x}\omega^*(n)^2$ with $x=10^k$ and $2\le k\le10$ using Mathematica. In view of the relation
\[M_2(x)=\int_{1}^{x}\frac{S_2(t)}{t}\,dt+O(1),\]
we also calculated the values of $S_2(x):=(1/x)\sum_{[p-1,q-1]\le x}1$ for $x$ in the same range. These values are recorded in the table below.

	\begin{table}[ht]
		\caption{Numerical values of $M_2(10^k)$ and $S_2(10^k)$}
		\label{Ta:xM_2(x)&xS_2(x)}
		\begin{tabular}{|l|l|l|l|l|} \hline
			~$k$  & $M_2(10^k)$ &$S_2(10^k)$&$3\log 10^k-6$&$3.2(1-1/\log 10^k)$\\ \hline
			~$2$ & 9.71&2.42&7.82 &2.51\\
			~$3$ & 15.530 &2.624&14.723 &2.737\\
			~$4$ & 21.9128&2.8175& 21.6310&2.8526\\
			~$5$ & 28.49311&2.88636& 28.53878&2.92205\\
			~$6$ & 35.261891&2.950910& 35.446532&2.968376\\
			~$7$ & 42.1296839&2.9923851& 42.3542870&3.0014654\\
			~$8$ & 49.07181351&3.02166709&49.26204223 &3.02628221 \\
			~$9$ & 56.067311859&3.043042188& 56.169797511&3.045584184\\
			~$10$ & 63.1033824202 &3.0595625181&63.0775527898 &3.0610257658\\
						\hline
		\end{tabular}
	\end{table}


The $M_2$ numbers in  Table \ref{Ta:xM_2(x)&xS_2(x)} seem to fit nicely with $3\log x-6$, and the $S_2$ numbers may fit with $3.2(1-1/\log x)$.  Perhaps $C\approx 3.1$?

\section{The level sets of $\omega^*(n)$}
\label{S:LS}
For $x,y\ge1$, let $N(x,y)\colonequals\#\{n \le x\colon\omega^*(n) \ge y\}$. The following theorem provides upper and lower bounds for $N(x,y)$.
\begin{theorem}
	\label{th:N(x,y)}
There exists a suitable constant $c>0$ such that 
\[\left\lfloor\frac{x}{y^{c\log\log y}}\right\rfloor\le N(x, y) \ll \frac{x\log y}{y}\]
for all $x\ge1$ and all sufficiently large $y$.
\end{theorem}
\begin{proof}
The lower bound follows immediately from \cite[Proposition 10]{APR}, which asserts that there exists some constant $c_0>0$ such that for all $z>100$, there is some positive integer $m_z<z$ with $\omega^*(m_z)>e^{c_0\log z/\log\log z}$. Taking $z=y^{c\log\log y}$ with some suitable constant $c>0$, we have $\omega^*(m_z)>y$ and hence
\[N(x,y)\ge\left\lfloor\frac{x}{m_z}\right\rfloor\ge\left\lfloor\frac{x}{y^{c\log\log y}}\right\rfloor.\]
\par To prove the upper bound, we first note that since the average of $\omega^*(n)$ for $n\le x$ is $\log\log x + O(1)$, it follows that $N(x,y)\ll x \log\log x /y$. So we have the desired upper bound when $y > (\log x)^{.05}$, say. Assume now that $y \le (\log x)^{.05}$, and let $z = \exp(y^{19})$, so that $z \le\exp((\log x)^{0.95})= x^{o(1)}$. There are two possibilities for $n$ counted by $N(x,y)$:
\begin{enumerate}
	\item $n$ is divisible by a shifted prime $p- 1>z$,
	\item $n$ is divisible by at least $y$ shifted primes $p-1\le z$.
\end{enumerate}
By \cite[Theorem 1.2]{MPP}, the count of the numbers in (1) is $\ll x/(\log z)^{\beta+o(1)}$, where $\beta\colonequals1-(1+\log\log 2)/\log 2$ is the Erd\H os--Ford--Tenenbaum constant.  Since $19\beta > 1$, the count in this case is $\ll x\log y/y$. For (2), let $\omega^*_z(n)$ denote the number of shifted primes $p-1\le z$ with $(p-1)\mid n$. It is easily seen that the average value of $\omega^*_z(n)$ for $n \le x$ is $\log\log z + O(1)$. Thus, the count in this case is $\ll x \log\log z / y\ll x\log y/y$. Adding up the bounds for the counts in both cases yields the desired upper bound for $N(x,y)$.
\end{proof}
\par Now we study the $k$-level set $\mathcal{L}_k\colonequals\{n\in\NN\colon \omega^*(n)=k\}$ for each $k\in\NN$. It is clear that 
\[N(x,y)=\sum_{k\ge y}\#(\mathcal{L}_k\cap[1,x]).\]
We shall show that each $\mathcal{L}_k$ has a positive natural density $\delta_k$, which is defined by
\begin{equation}
	\label{def:delta_k}
\delta_k\colonequals\lim_{x\to\infty}\frac{\#(\mathcal{L}_k\cap[1,x])}{x}.
\end{equation}

\begin{theorem}
	\label{th:delta_k}
For every $k\in\NN$, the $k$-level set  $\mathcal{L}_k$ admits a positive natural density $\delta_k$. Moreover, we have $\sum_k\delta_k=1$.
\end{theorem}

We first show that each $\mathcal{L}_k$ is nonempty.

\begin{lemma}
	\label{lem:L_k}
For every $k\in\NN$, we have $\mathcal{L}_k\ne\emptyset$.
\end{lemma}

\begin{proof}
Note that $\mathcal{L}_1=\NN\setminus2\NN$ and $2\in\mathcal{L}_2$. So we may assume that $k\ge2$, so that $\mathcal{L}_k\subseteq2\NN$. We shall show that for any $n\in2\NN$, there exists an integral multiple $m\in\NN$ of $n$ such that $\omega^*(m)=\omega^*(n)+1$. The lemma would then follow from this result in an inductive manner.
\par To prove this, we fix $n\in2\NN$ and consider
\[\mathcal{P}_2(x)\colonequals\{2<p\le x\colon\Omega((p-1)/2)\le 2\text{~and~}P^-((p-1)/2)>x^{3/11}\}.\]
(The notation here is standard, signifying that $(p-1)/2$ is either prime or the product of two primes,
and this prime or primes are $>x^{3/11}$.)
By \cite[Theorem 25.11]{FI}, we have $\#\mathcal{P}_2(x)\gg x/(\log x)^2$ for all sufficiently large $x$. We wish to find some large $p\in\mathcal{P}_2(x)$ with $\omega^*(n(p-1)/2)=\omega^*(n)+1$. To this end, we shall show that the number of those $p\in\mathcal{P}_2(x)$ which do not have this property is $O(x\log\log x/(\log x)^3)$. Note that if $p\in\mathcal{P}_2(x)$ does not possess this property, then we can find $a\mid n$ and $b\mid (p-1)/2$ with $a,b>1$ such that $ab+1$ is a
prime not equal to $p$.

There are two possibilities: (i) $b=(p-1)/2$ and $ab+1$ is prime with $a>2$
and (ii) $p-1=2qr$ with $q,r$ primes in $(x^{3/11},x^{8/11}/2)$ and $aq+1$ is prime.

Case (i) is simple.  Fix $a\mid n$ with $a>2$.
The number of integers $b\le x$ with $P^-(b)>x^{3/11}$ and
both $2b+1$ and $ab+1$ are prime is $\ll x/(\log x)^3$.  (The implied constant
depends on $a$ but there is a bounded number of choices for $a$.)

Now we consider Case (ii). Again, let us fix $a\mid n$. For any prime $q\in(x^{3/11},x^{8/11}/2)$, the number of primes $b< x/2q$ such that both $ab+1$ and $2qb+1$ are prime is
\[\ll\frac{x}{q(\log x)^3}\prod_{r\mid (2q-a)}\left(1-\frac{1}{r}\right)^{-1}\ll\frac{\log\log q}{q}\cdot\frac{x}{(\log x)^3}.\]
Summing this bound over all $q\in(x^{3/11},x^{8/11}/2)$ and $a\mid n$, we see that the number of choices of $p$ with $b$ in Case (ii) is $\ll x\log\log x/(\log x)^3$.
This completes the proof.
\end{proof}

\par We are now ready to prove Theorem \ref{th:delta_k}.

\begin{proof}[Proof of Theorem \ref{th:delta_k}]
The case $k=1$ is obvious, since the level set $\mathcal L_1$
consists of the odd numbers, so that $\delta_1=1/2$. Now let us fix $k\ge2$. Then $\mathcal{L}_k\subseteq2\NN$. We define an equivalence relation $\simeq$ on $\mathcal{L}_k$ by declaring that $m\simeq n$ if and only if $m$ and $n$ have exactly the same set of shifted prime divisors. Let $\mathcal{C}_k$ be the set of all equivalence classes\footnote{~In \cite{PW}, the class containing $n$ is denoted $\mathcal S_n$.} $\langle n\rangle$ of $\mathcal{L}_k$ under $\simeq$. Then 
\begin{equation}
	\label{eq:L_k}
\mathcal{L}_k=\bigcup_{\langle n\rangle\in\mathcal{C}_k}\langle n\rangle.
\end{equation} 
It is known \cite[Theorem 3]{EW} that each $\langle n\rangle$ has a positive natural density. Thus, if natural density were countably additive, then we would be able to conclude that $\delta_k$ exists and equals the sum of the natural densities of the sets $\langle n\rangle\in\mathcal{C}_k$. Since Lemma \ref{lem:L_k} implies that $\mathcal{C}_k\ne\emptyset$, we would also have $\delta_k>0$. Unfortunately, $\#\mathcal{C}_k$ may be infinite and natural density is only finitely additive. 

To overcome this issue we appeal to the following elementary result.
\begin{lemma}
\label{lem:dense}
Let $\mathcal A_1,\mathcal A_2,\dots$ be an infinite sequence of pairwise disjoint subsets of $\NN$, such 
that each $\mathcal A_i$
has a natural density $\delta(\mathcal A_i)$.  If the upper asymptotic density of
$\bigcup_{i>j}\mathcal A_i$ tends to $0$ as $j\to\infty$, then the density of $\bigcup_{i\ge1}\mathcal A_i$
exists and
\[
\delta\Big(\bigcup_{i\ge1}\mathcal A_i\Big)=\sum_{i\ge 1}\delta(\mathcal A_i).
\]
\end{lemma}
This result can be applied to the sets $\langle n\rangle\in\mathcal C_k$, since if there are infinitely
many, then for any fixed $y$ all but finitely many have $n$ divisible by a shifted prime $p-1>y$.
 Appealing to \cite[Theorem 2]{EW}, the union of these sets has
upper density tending to 0 as $y\to\infty$.  Thus, to complete the proof, we now have
\[
\delta_k=\sum_{\langle n\rangle\in\mathcal C_k}\delta(\langle n\rangle)
~\hbox{ and }~\sum_k\delta_k=1.
\]
\end{proof}

Here are some exact counts of the level sets $\mathcal L_k$ for $k\le 11$.

\begin{table}[ht]
		\caption{Exact counts of level sets for $k<12$}
		\label{Ta:delta}
		\centering
		\begin{tabular}{|r|r|r|r|r||r|}\hline
		&&&&&\\
			~$k$ &$10^4$ & $10^6$ &$10^8$ & $10^{10}$ & $\approx\delta_k$\\\hline
			~1&$5{,}000$ & $500{,}000$ & $50{,}000{,}000$ & $5{,}000{,}000{,}000$ &  .5\\ 
			~$2$ &834 & 77{,}696 & 7{,}436{,}825 & 720{,}726{,}912 &.070\\
			~$3$& 965 & 91{,}602 & 8{,}826{,}498 & 859{,}002{,}140 &.084\\
			~$4$& 877 &79{,}986 & 7{,}691{,}971 &  748{,}412{,}490&.074\\
			~$5$& 612 & 59{,}518 & 5{,}684{,}323 &555{,}900{,}984  &.055\\
			~$6$ & 456 & 40{,}641 & 4{,}031{,}009 & 401{,}146{,}301 &.040\\
			~$7$ & 287 & 29{,}565 &3{,}016{,}881 & 300{,}330{,}932 &.030\\
			~$8$ & 202 & 23{,}190 & 2{,}324{,}769 & 233{,}611{,}502 &.023 \\
			~$9$ & 153 & 17{,}914 & 1{,}800{,}298 & 182{,}793{,}491 &.018\\
			~$10$ & 159 & 13{,}899 & 1{,}401{,}307 &  144{,}740{,}573&.015\\
			~11 & 103& 10{,}487 & 1{,}131{,}836 &  118{,}302{,}267&.012\\
			$\ge12$& 352& 55{,}682&6{,}654{,}283&735{,}032{,}408&\\
						\hline
		\end{tabular}
	\end{table}

The largest values of $k$ encountered here up to the various bounds:
$10^4$:  28, $10^6$:  86, $10^8$: 247, $10^{10}$: 618.

It is curious that $\delta_2$ is apparently smaller than $\delta_3$.  This is partially explained
by the fact that $\mathcal C_2$ is the single equivalence class $\langle 2\rangle$, while
$\mathcal C_3$ contains $\langle4\rangle,\langle6\rangle,\langle10\rangle$ and infinitely
many other classes (see \cite[Theorem 3]{PW}).  While $\delta(\langle2\rangle)$ is the
largest $\delta(\langle n\rangle)$ for $n$ even (proved in \cite{PW}, slightly improving
an earlier result of Sunseri), the smaller ones apparently unite to
surpass the single larger density.
Perhaps though the densities $\delta_k$ are monotone for $k\ge3$.  
 It would be good to have some sort of asymptotic inequalities for these densities, and a result in this direction is produced in the next section.

\section{A lower bound for $\delta(\langle n\rangle)$}

Let $n\in2\NN$ and consider the equivalence class $\langle n\rangle$ of $\NN$ under the same relation $\simeq$ as introduced in the proof of Theorem \ref{th:delta_k} above. Suppose that $n=\min\langle n\rangle$. In other words, $n$ is the least common multiple of all shifted prime divisors of $n$. We clearly have $\delta(\langle n\rangle)<1/n$. Erd\H{o}s and Wagstaff \cite{EW} asked what a positive lower bound could be for $\delta(\langle n\rangle)$. The following theorem provides such a lower bound.

\begin{theorem}
	\label{th:delta([n])}
Let $n\in2\NN$ be such that $n=\min\langle n\rangle$. Then 
\[\delta(\langle n\rangle)\ge\frac{1}{n^{O(\tau(n))}}.\]
\end{theorem}
\begin{proof}
We follow the proof of \cite[Theorem 3]{EW} on the existence and positivity of $\delta(\langle n\rangle)$. For any $a_1,...,a_r\in\NN$, denote by $T_n(a_1,...,a_r)$ the natural density of the set of multiples of $n$ which are not divisible by any $a_i$ for $1\le i\le r$. Explicitly, we have
\[T_n(a_1,...,a_r)=\sum_{j=0}^{r}(-1)^j\sum_{1\le i_1<\cdots<i_j\le r}\frac{1}{[n,a_{i_1},...,a_{i_j}]}.\]
By \cite[Eq. (2), p. 110]{EW}, we have
\[\frac{1}{n}T_n(a_1,...,a_{r+s})\ge T_n(a_1,...,a_r) T_n(a_{r+1},...,a_{r+s})\]
for any integers $r,s\ge0$ and any $a_1,...,a_{r+s}\in\NN$. From this inequality with $s=1$ it follows immediately by induction that
\begin{equation}
\label{eq:T_n}
T_n(a_1,...,a_r)\ge\frac{1}{n}\prod_{i=1}^{r}\left(1-\frac{n}{[n,a_i]}\right)^{1/m_i},
\end{equation}
where $m_i\colonequals\#\{1\le j\le r\colon a_j=a_i\}$. \par It suffices to prove the theorem for large values of $n$. Let $y\ge1$ be a parameter depending on $n$. Let $\mathscr{A}_1\colonequals \{[p-1,n]\colon p-1\le y\text{~and~}(p-1)\nmid n\}$ and $\mathscr{A}_2\colonequals \{[p-1,n]\colon p-1> y\text{~and~}(p-1)\nmid n\}$. Denote by  $\mathscr{B}(\mathscr{A}_2)$ the set of multiples of elements of $\mathscr{A}_2$. We arrange the elements of $\mathscr{A}_1$ as a strictly increasing sequence  $\{a_i\}_{i=1}^r$. The proof of \cite[Theorem 3]{EW} shows that 
\[\frac{1}{n}\delta(\langle n\rangle)\ge T_n(a_1,...,a_r)\left(\frac{1}{n}-\delta(\mathscr{B}(\mathscr{A}_2))\right)>0,\]
provided that $y$ is sufficiently large in terms of $n$. To get a positive lower bound for $\delta(\langle n\rangle)$, it suffices to obtain a positive lower bound for $T_n(a_1,...,a_r)$ and an upper bound $<1/n$ for $\delta(\mathscr{B}(\mathscr{A}_2))$. Take $y=e^{n^{1/\beta}}$.
By \cite[Theorem 1.2]{MPP} we have
\[\delta(\mathscr{B}(\mathscr{A}_2))\ll\frac{1}{(\log y)^{\beta}\sqrt{\log\log y}}\ll\frac{1}{n\sqrt{\log n}}.\]
 To handle $T_n(a_1,...,a_r)$, we appeal to \eqref{eq:T_n} to obtain
\begin{align*}
T_n(a_1,...,a_r)&\ge\frac{1}{n}\prod_{d\mid n}\prod_{\substack{p\le y+1\\(p-1,n)=d\\(p-1)\nmid n}}\left(1-\frac{d}{p-1}\right)\\
&\ge\frac1n\exp\Bigg(-\sum_{d\,|\,n}\sum_{\substack{p\le y+1\\(p-1,n)=d\\(p-1)\,\nmid\, n}}\frac d{p-1}\Bigg)\\
&=\frac{e^{\omega^*(n)}}n\exp\Bigg(-\sum_{d\,|\,n}\sum_{\substack{p\le y+1\\(p-1,n)=d}}\frac d{p-1}\Bigg).
\end{align*}
If we replace $d/(p-1)$ with $d/p$, the error created in the double sum is $\ll\sigma(n)/n$,
where $\sigma$ is the sum-of-divisors function.  Thus,
\begin{equation}
\label{eq:doubsum}
T_n(a_1,...,a_r)\ge\frac{e^{\omega^*(n)}}n\exp\Bigg(-\sum_{d\,|\,n}\sum_{\substack{p\le y+1\\(p-1,n)=d}}\frac dp+O(\log\log n)\Bigg).
\end{equation}
\begin{lemma}
\label{lem:SW}
For each number $A>0$ there is a positive constant $\kappa$ such that for all
large $x$ and $d<(\log x)^A$, we have
\[
\sum_{\substack{d<p\le x\\p\,\equiv\, a\kern-5pt\pmod{d}}}\frac{1}{p}=\frac{\log\log x}{\varphi(d)}+E(d)+O(\exp(-\kappa(\log x)^{1/2})),
\]
for all $a$ coprime to $d$.  The number $E(d)$ satisfies $|E(d)|\ll\log(2d)/\varphi(d)$.
\end{lemma}
This follows from the Siegel--Walfisz theorem, where the estimation for $E(d)$ appears
in works of Norton and Pomerance, see  \cite[Lemma 2.1]{MM}.

 Consider the double sum in \eqref{eq:doubsum}. It follows that
\begin{align*}
\sum_{d\mid n}&d\sum_{\substack{p\le y+1\\(p-1,n)=d}}\frac{1}{p}=\sum_{cd\mid n}\mu(c)d\sum_{\substack{p\le y+1\\p\equiv 1\Mod{cd}}}\frac{1}{p}\\
&=\sum_{cd\mid n}\mu(c)d\Bigg(\frac{\log\log y}{\varphi(cd)}
+E(cd)+O(\exp(-\kappa(\log y)^{1/2}))\Bigg)\\
&=\tau(n)\log\log y+\sum_{m\,\mid\, n}\varphi(m)E(m)+O(n^2\exp(-\kappa(\log y)^{1/2}))\\
&\le \tau(n)(\log\log y+O(\log n))=O(\tau(n)\log n).
\end{align*}
Hence, we have
\[T_n(a_1,...,a_r)\ge\frac{e^{\omega^*(n)}}{n^{O(\tau(n))}}=\frac1{n^{O(\tau(n))}},\]
the last estimate coming from $\omega^*(n)\le\tau(n)$.
Combining the above estimate with that for $\delta(\mathscr{B}(\mathscr{A}_2))$
completes the proof.
\end{proof}

\section{Higher moments of $\omega^*(n)$}
\label{S:3M}
For every $k\in\NN$, we define the $k$th moment of $\omega^*(n)$ by
\[M_{k}(x)\colonequals\frac{1}{x}\sum_{n\le x}\omega^*(n)^k.\]
By opening the power and reversing the order of summation we have
\begin{equation}
	\label{eq:M_k}
M_k(x)=\frac{1}{x}\sum_{[p_1-1,...,p_k-1]\le x}\left\lfloor\frac{x}{[p_1-1,...,p_k-1]}\right\rfloor.
\end{equation}
This shows that $M_k(x)$ is closely related to 
\[S_{k}(x)\colonequals\frac{1}{x}\sum_{[p_1-1,...,p_k-1]\le x}1.\]
In fact, if $S_k(x)\ll (\log x)^{c_k}$, then a partial summation argument applied to the
upper bound in \eqref{eq:M_k} afforded by removing the floor function shows that
$M_k(x)\ll (\log x)^{c_k+1}$.  A similar argument shows that a lower bound for $S_k(x)$
implies one for $M_k(x)$.

For $k\ge2$, it is natural to relate the function $\omega^*(n)^k$ to $\tau(n)^k$. It is well-known that for every $k\ge1$, one has
\[\frac{1}{x}\sum_{n\le x}\tau(n)^k\sim\frac{1}{(2^k-1)!}\prod_{p}\left(\left(1-\frac{1}{p}\right)^{2^k}\sum_{\nu\ge0}\frac{(\nu+1)^k}{p^{\nu}}\right)(\log x)^{2^k-1},\]
in contrast to 
\[\frac{1}{x}\sum_{n\le x}\omega(n)^k\sim(\log\log x)^k.\]
Comparing $\omega^*$ with $\tau$ and taking the primality conditions into account, one may conjecture that $M_k(x)\sim \mu_k(\log x)^{2^k-k-1}$ for every $k\ge2$, where $\mu_k>0$ is a constant depending on $k$. Similarly, one may also conjecture that  $S_k(x)\sim (2^k-k-1)\mu_k(\log x)^{2^k-k-2}$ for every $k\ge2$ with the same constant $\mu_k$. As in the case $k=2$, we have the upper and lower bounds for $M_3(x)$  of the conjectured magnitude.\begin{theorem}
\label{th:3M}
We have $M_3(x)\asymp(\log x)^4$ for all $x\ge2$.
\end{theorem}
The upper and lower bounds will be proved by using different types of arguments. The rest of this section will be devoted to proving the upper bound $M_3(x)\ll(\log x)^4$, with the proof of the lower bound $M_3(x)\gg(\log x)^4$ given in Section \ref{sec:low}.

We begin with some lemmas.  The first is a variant of \cite[Lemma 2.7]{MM}.
\begin{lemma}
\label{lem:abef}
Uniformly for coprime integers $e,f$ in $[1,x]$,
\begin{equation}
\label{eq:abef}
\sum_{\substack{a,b\le x\\(ae,bf)=1\\ae\ne bf}}\frac1{ab}\prod_{p\,|\,ab(ae-bf)}\left(1+\frac1p\right)
\ll (\log x)^2.
\end{equation}
\end{lemma}
\begin{proof}
First note that the product contributes at most a factor of magnitude $\log\log x$ to the sum,
so the result holds trivially if either $a$ or $b$ is bounded by $x^{1/\log\log x}$.  Hence,
we may assume that $a,b>x^{1/\log\log x}$.  Further, every integer $n\le x$ has
$<\log x$ prime divisors, so that
\[
\prod_{\substack{p\,|\,n\\p>(\log x)^{1/2}}}\left(1+\frac1p\right)\ll 1,
\]
uniformly.   Let $u$ be the product of all primes $p\le(\log x)^{1/2}$.
Thus, we may restrict the primes $p$ in the product in the lemma to those that also divide $u$.
We have the expression in \eqref{eq:abef} is
\begin{equation}
\label{eq:abef2}
\ll\sum_{\substack{x^{1/\log\log x}<a,b\le x\\(ae,bf)=1\\ae\ne bf}}\frac1{ab}\sum_{\substack{j\,|\,u\\j\,|\,ab(ae-bf)}}\frac1j
\le\sum_{j\,|\,u}\frac1j\sum_{\substack{j<a,b\le x\\j\,|\,ae-bf\\(ae,bf)=1\\ae\ne bf}}\frac1{ab}.
\end{equation}
(Note that we assume here that $a,b >j$, since they are $>x^{1/\log\log x}$ and
$j\le u\le\exp((1+o(1))(\log x)^{1/2})$.)
For $p\mid j$ with $j\mid ab(ae-bf)$,
we have either $a\equiv0\pmod p$, $b\equiv0\pmod p$, or $a\equiv bfe^{-1}\pmod p$
(if $p\mid e$ then $p\nmid ae-bf$).  Since $j$ is squarefree, there are at most
$3^{\omega(j)}j$ pairs $a,b\pmod j$ with $j\mid ab(ae-bf)$.  For a fixed pair of residues
(mod~$j$) that we have here, the sum of $1/ab$ in this class is $\ll(\log x)^2/j^2$ uniformly, so the
total contribution in the last sum in \eqref{eq:abef2} is $\ll 3^{\omega(j)}(\log x)^2/j$.  We thus have the last double sum in
\eqref{eq:abef2} is
\[
\ll\sum_{j\mid u}\frac{3^{\omega(j)}(\log x)^2}{j^2}=(\log x)^2\prod_{p\le(\log x)^{1/2}}\left(1+\frac3{p^2}\right)\ll(\log x)^2,
\]
which completes the proof of the lemma.
\end{proof}
 
\begin{lemma}
\label{lem:ulem}
Uniformly for $1\le u<x$ we have
\[
\sum_{\substack{q\le x\\q\,\equiv\,1\kern-7pt\pmod u}}\frac{\tau((q-1)/u)}{\varphi((q-1)/u)}\ll\frac u{\varphi(u)}\log x.
\]
\end{lemma}
\begin{proof}
The result holds trivially when $x\le 2$ or $u\ge x/2$, so assume that $x>2$ and $u<x/2$.
(In fact, the lemma follows from a trivial argument if $u>x/\exp((\log x)^{1/2})$, but we
won't use this.)  We first consider
\[
T=\sum_{\substack{q\le x\\q\,\equiv\,1\kern-7pt\pmod u}}\frac{\tau((q-1)/u)(q-1)/u}{\varphi((q-1)/u)}.
\]
Using that
\[
\frac n{\varphi(n)}=\sum_{d\mid n}\frac{\mu^2(d)}{\varphi(d)},
\]
we have
\[
T\le\sum_{d<x}\frac1{\varphi(d)}\sum_{\substack{q\le x\\q\,\equiv\,1\kern-7pt\pmod{du}}}
\tau((q-1)/u).
\]
Using the maximal order of the divisor function and that $q$ is an integer that is 1 (mod $ud$) and
$>ud$, the contribution to $T$ from a particular
number $d$ is $\ll (x/du)(x/u)^\epsilon$, so the contribution to $T$ from numbers $d>(x/u)^{1/4}$
is $\ll (x/u)^{3/4+\epsilon}<(x/u)^{4/5}$, say.  We also use that $\tau((q-1)/u)$ is at most twice the number of
divisors $j\mid(q-1)/u$ with $j\le(x/u)^{1/2}$.  Thus,
\[
T\ll\sum_{d\le(x/u)^{1/4}}\frac1{\varphi(d)}\sum_{j\le(x/u)^{1/2}}\sum_{\substack{q\le x\\q\,\equiv\,1\kern-7pt\pmod{[j,d]u}}}1+(x/u)^{4/5}.
\]
We have $[j,d]\le(x/u)^{3/4}$ and so $x/([j,d]u)\ge(x/u)^{1/4}$ and the inner sum here is $\ll x/(\varphi([j,d]u)\log(x/u))$.
Now $[j,d]=jd/i$, where $i=(j,d)$, so
\[
T\ll\sum_{d\le(x/u)^{1/4}}\frac1{\varphi(d)}\sum_{i\mid d}\sum_{k\le(x/u)^{1/2}/i}\frac x{\varphi(d)\varphi(u)\varphi(k)\log(x/u)}+(x/u)^{4/5}.
\]
The sum of $1/\varphi(k)$ in the indicated range is $\ll\log(x/u)$, so
\[
T\ll\sum_{d\le (x/u)^{1/4}}\frac{x\tau(d)}{\varphi(u)\varphi(d)^2}+(x/u)^{4/5}\ll\frac x{\varphi(u)}.
	\]
	It immediately follows that
\[
\sum_{\substack{q\le x\\q\,\equiv\,1\kern-7pt\pmod u}}\frac{\tau((q-1)/u)(q-1)}{\varphi((q-1)/u)}
	\ll\frac{ux}{\varphi(u)}
	\]
and the lemma follows by partial summation.
\end{proof}

At this point we find it convenient to reprise the upper bound proof for $k=2$ from \cite{MM} since we take
a slightly different perspective, the proof is short, and the case $k=3$ follows with similar tools.
We show that
\begin{equation}
\label{eq:k=2}
S_2(x)\ll 1.
\end{equation}
We are to count pairs of primes $p,q$ with $[p-1,q-1]\le x$.  Let 
\[
d=(p-1,q-1),~~p-1=ad,~~q-1=bd.
\]
The case $p=q$ has the count $O(x/\log x)$, so we may assume that $a\ne b$.
So, we are counting triples $a,b,d$ with $(a,b)=1$, $a\ne b$, $abd\le x$, with $ad+1,bd+1$ both prime.
First suppose that $d=\max\{a,b,d\}$.   Since $abd\le x$, 
we have $ab\le x^{2/3}$.  For a given choice of $a,b$, the number of choices for $d\le x/ab$
with $ad+1,bd+1$ both prime is
\[
\ll\frac x{ab(\log x)^2}\prod_{\ell\,|\,ab(a-b)}\left(1+\frac1\ell\right),
\]
where $\ell$ runs over primes.  This follows from the upper bound in either Brun's or
Selberg's sieve.  Lemma \ref{lem:abef} in the case $e=f=1$ completes the proof of
\eqref{eq:k=2} in this case.

Now assume that $a=\max\{a,b,d\}$.  Then $bd=q-1\le x^{2/3}$.  For a given prime
$q\le x^{2/3}+1$ and a divisor $d$ of $q-1$, we count values of $a\le x/(q-1)$ with
$ad+1$ prime.  By the Brun--Titchmarsh inequality, the number of such values of $a$
is 
\[
\ll\frac d{\varphi(d)}\frac x{(q-1)\log x}\le \frac x{\varphi(q-1)\log x}. 
\]
So, in all, there are $\ll \tau(q-1)x/(\varphi(q-1)\log x)$ choices for $a$.
 Lemma \ref{lem:ulem} in
the case $u=1$ completes the proof of \eqref{eq:k=2} in this case.  The last case
$b=\max\{a,b,d\}$ is completely symmetric with the case just considered, so we
are done.

Now we prove the upper bound $M_3(x)\ll(\log x)^4$ asserted in Theorem \ref{th:3M}.  For this it is sufficient to prove that
\begin{equation}
	\label{eq:pqr}
	S_3(x)=\frac{1}{x}\sum_{[p-1,\,q-1,\,r-1]\,\le\, x}1\ll (\log x)^3,
\end{equation}
where $p,q,r$ run over prime numbers.  From the case of the second moment, we may
assume that $p,q,r$ are distinct.  Note that
\[
[p-1,q-1,r-1]=abcde\kern-1pt f\kern-1pt g,
\]
where
\[
g=\gcd(p-1,q-1,r-1),
\]\[
dg=\gcd(p-1,q-1),~~eg=\gcd(p-1,r-1),~~f\kern-1pt g=\gcd(q-1,r-1),
\]\[
a=(p-1)/deg,~~b=(q-1)/df\kern-1pt g,~~c=(r-1)/e\kern-1pt f\kern-2pt g.
\]
Note that we have $a,b,c$ pairwise coprime, as well as $d,e,f$.  Also,
\[
\gcd(ae,b\kern-1pt f)=1,~~\gcd(ad,ce)=1,~~\gcd(bd,c\kern-1pt f)=1.
\]

To prove \eqref{eq:pqr}, we consider 7 cases depending on the largest of
$a,\dots,g$.  By symmetry these collapse to 3 cases:
\[
\max\{a,\dots,g\}=a,d,\,\hbox{or }g.
\]
Beginning with the max being $a$, first choose a prime $r$ and a factorization of
$r-1$ as $ce\kern-1pt f\kern-1pt g$.  Next choose a prime $q$ with $q\equiv1\pmod{f\kern-1pt g}$ and take
a factorization of $(q-1)/f\kern-1pt g$ as $bd$.  Finally, let $a\le x/bcde\kern-1pt f\kern-1pt g$ with $adeg+1$ prime.
The number of choices for $a$ is 
\[
\ll\frac {x}{bcde\kern-1pt f\kern-1pt g\log x}\frac{deg}{\varphi(deg)}.
\]
The number of choices for $c,e,f,g$ is $\tau_4(r-1)$.  Given an ordered factorization
$ce\kern-1pt f\kern-1pt g$ of $r-1$, let $u=u_{r{-}1}=f\kern-1pt g$.  The
number of choices for $b,d$ is $\tau((q-1)/u)$.  Thus, the total number of choices in this case is
\[
\ll \sum_{r<x}\frac{\tau_4(r-1)}{\varphi(r{-}1)}\sum_{\substack{q<x\\q\,\equiv\,1\kern-7pt\pmod{u_{r-1}}}}
\frac{\tau((q-1)/u)}{\varphi((q-1)/u)}\frac x{\log x}.
\]
Using Lemma \ref{lem:ulem}, we have the number of choices
\begin{equation}
\label{eq:pp}
\ll x\sum_{r<x}\frac{\tau_4(r-1)}{\varphi(r-1)}\frac{u}{\varphi(u)}
\le x\sum_{r<x}\frac{\tau_4(r-1)(r-1)}{\varphi(r-1)^2}.
\end{equation}
We now appeal to \cite[Theorem 1.2]{PP16} or \cite[Corollary 1.2]{PP20} from which we see this
last sum is $\ll(\log x)^3$.  This completes the proof when
$a=\max\{a,\dots,g\}$.  

Now assume that $d=\max\{a,\dots,g\}$.  We choose a prime $r<x$ and a
factorization $ce\kern-1pt f\kern-1pt g$ of $r-1$.  We then choose $a,b$ with $ab(r-1)\le x^{6/7}$.
We now let $d$ run up to $x/(ab(r-1))$ with $adeg+1$ and $bdf\kern-1pt g+1$ prime.
The number of choices is
\[
\ll\frac{x\tau_4(r-1)}{ab(r-1)(\log x)^2}\prod_{\ell\mid abe\kern-1pt f(bf-ae)}\Big(1+\frac1\ell\Big)
\prod_{\ell\mid g}\Big(1+\frac1\ell\Big),
\]
where $\ell$ runs over prime numbers.  We can absorb the part of the product coming
from $\ell\mid e\kern-1pt f$ and $\ell\mid g$ into the main term, getting
\[
\frac{x\tau_4(r-1)(r-1)}{\varphi(r-1)^2ab(\log x)^2}\prod_{\ell\mid ab(bf-ae)}\Big(1+\frac1\ell\Big).
\]
Note that this final product is finite, since $(ae,bf)=1$ and $ae\ne bf$. (If $a=b=e=f=1$, then
one has $p=q$, a possibility we ruled out).  
Lemma \ref{lem:abef}  and then the argument as in \eqref{eq:pp} completes the proof of the case when $d$ is the maximum of $a,\dots,g$.

We now consider the case that $g=\max\{a,\dots,g\}$, 
which is quite similar to the previous case.
For a given choice of $a,\dots,f$, we have $a\dots f\le x^{6/7}$, 
so the number of values of $g\le x/a\dots f$ with
$adeg+1,bdf\kern-1pt g+1,ce\kern-1pt f\kern-1pt g+1$ all prime is
\[
\ll \frac x{A(\log x)^3}\prod_{\ell\mid AE}\Big(1+\frac2\ell\Big),
\]
where
\[
A=abcde\kern-1pt f,\quad E=(ae-bf)(ad-cf)(bd-ce).
\]
Without the product, the sum of $1/A$ is $O((\log x)^6)$.  We would like to
show the same estimate holds with the product included.  Note however that the
product is in the worst case $O((\log\log x)^2)$, so our result holds trivially if any
of $a,\dots,f$ is $\le x^{1/(\log\log x)^2}$.  We thus assume they are all $>x^{1/(\log\log x)^2}$.
Further, as in the proof of Lemma \ref{lem:abef}, let $u$ be the product of all primes $\ell\le(\log x)^{1/2}$.
We may restrict primes $\ell$ in the product to  such primes.  We wish to estimate
\[
\sum_{j\mid u}\frac{2^{\omega(j)}}j\sum_{\substack{a,\dots,f<x\\j\mid AE}}\frac1{a\dots f}.
\]
Note that $AE$ is the product of 9 expressions, so that in the inner sum,
the 6-tuple $(a,\dots,f)$ lies in $\le 9^{\omega(j)}j^5$ residue classes mod $j$.
For each one of these classes the inner sum is $\ll(\log x)^6/j^6$, so we have
\[
\sum_{j\mid u}\frac{2^{\omega(j)}}j\frac{9^{\omega(j)}}{j}(\log x)^6=(\log x)^6\sum_{j\mid u}\frac{18^{\omega(j)}}{j^2}\ll(\log x)^6.
\]
This completes our proof of the upper bound in Theorem \ref{th:3M}.

\section{A lower bound for the third moment}
\label{sec:low}

By \eqref{eq:M_k}, we have
\[
M_3(x)\ge\frac12\sum_{[p-1,q-1,r-1]\le x}\frac1{[p-1,q-1,r-1]}.
\]
Thus, we wish to show that
\begin{equation}
\label{eq:3low}
\sum_{[p-1,q-1,r-1]\le x}\frac1{[p-1,q-1,r-1]}\gg(\log x)^4.
\end{equation}
We restrict to the case that $p,q,r$ are distinct primes, noting that
the complementary case is negligible. We use the identity
\[
\frac1{[p-1,q-1,r-1]}=\sum_{\substack{u\,|\,r-1\\u\,|\,[p-1,q-1]}}\frac{\varphi(u)}{[p-1,q-1](r-1)}.
\]
Let
\begin{equation}
\label{eq:M2u}
M_2(x;u):=\sum_{\substack{[p-1,q-1]\le x\\u\,|\,[p-1,q-1]}}\frac1{[p-1,q-1]}.
\end{equation}
We thus have that
\begin{equation}
\label{eq:m3}
M_3(x)\ge\frac12\sum_{u\le x^{1/3}}\sum_{\substack{r\le x^{1/3}\\u\,|\,r-1}}\frac1{r-1}M_2(x^{2/3};u).
\end{equation}
Our goal then is to obtain a lower bound for $M_2(x^{2/3},u)$ and use that in \eqref{eq:m3}.

Helpful will be a tool from \cite{AGP}, namely Theorem 2.1:  {\it For each $\varepsilon>0$ there are numbers $\delta>0$ and $x_0$, such that if
$x>x_0$, $k<x^\delta$, and $(a,k)=1$, then
\begin{equation}
\label{eq:AGP}
\Bigg|\sum_{\substack{p\le y\\p\,\equiv\,a\kern-5pt\pmod k}}\log p-\frac y{\varphi(k)}\Bigg|
\le\varepsilon\frac y{\varphi(k)},
\end{equation}
for all $y\ge x$, except possibly for those $k$ divisible by a certain number
$k_0(x)>\log x$.}

If $k_0(x)$ should exist and is divisible by a prime $>(1/3)\log\log x$, let $s=s(x)$ be the
largest such prime.  Otherwise, let $s=s(x)$ be the least prime $>(1/3)\log\log x$.
Note that if $x$ is sufficiently large and $k_0(x)$ exists and is not divisible by any prime $>(1/3)\log\log x$,
then $k_0(x)$ must be divisible by the cube of some prime.  Indeed
\[
\prod_{\ell\le(1/3)\log\log x}\ell^2=(\log x)^{2/3+o(1)},
\]
so this product is smaller than $k_0(x)$.  In particular, if $x$ is large
and $k$ is cube-free, then \eqref{eq:AGP} holds whenever $s\nmid k$.
\begin{corollary}
\label{cor:apr}
Suppose that $\varepsilon=1/4$ and we have the corresponding number $\delta$ as
above.  If $x$ is sufficiently large, $k<x^\delta$ is cube-free and not divisible by $s(x)$, then
\[
\sum_{\substack{x<p\le x^2\\p\,\equiv\,a\kern-5pt\pmod k}}\frac{1}p\gg\frac1{\varphi(k)},
\]
uniformly.
\end{corollary}
This follows instantly from the above theorem, namely \cite[Theorem 2.1]{AGP}, by either
partial summation or a dyadic summation.

Let $\varphi_2(n)$ be the multiplicative function with value at a prime power $\ell^j$ equal
to $\ell^j(1-2/\ell)$.
\begin{proposition}
\label{prop:m2u}
Let $\delta$ be the corresponding constant for $\varepsilon=1/4$.
Suppose that $x$ is large and $u<x^{\delta/12}$ is squarefree and not divisible by $s=s(x^{1/6})$.
Then 
\[
M_2(x^{2/3};u)\gg\log x\sum_{\substack{u=u_1u_2u_3\\u_1u_2\,{\text odd}}}
\frac{u_3\varphi_2(u_1u_2)}{\varphi(u)^2},
\]
uniformly.
\end{proposition}
\begin{proof}
We have
\begin{equation}
\label{eq:m2/3}
M_2(x^{2/3};u)\ge\sum_{\substack{u=u_1u_2u_3\\u_1u_2\,{\rm odd}}}\sum_{\substack{p\le x^{1/3}
\\p\,\equiv\,1\kern-5pt\pmod{u_1u_3}\\(p-1,u_2s)=1}}
\sum_{\substack{q\le x^{1/3}\\
q\,\equiv\,1\kern-5pt\pmod{u_2u_3}\\(q-1,u_1s)=1}}\frac1{[p-1,q-1]}.
\end{equation}
We write $1/[p-1,q-1]$ as $(p-1,q-1)/(p-1)(q-1)$ and note that $u_3\mid(p-1,q-1)$.
Thus,
\[
\frac1{[p-1,q-1]}=u_3\sum_{d\,|\,(p-1,q-1)/u_3}\frac{\varphi(d)}{(p-1)(q-1)}.
\]
Thus, for a given choice of $u_1,u_2,u_3$, the contribution in \eqref{eq:m2/3} is at least
\begin{equation}
\label{eq:u1u2u3}
\sum_{\substack{d<x^{1/12}\\d\,{\rm squarefree}\\(d,s)=1}}u_3\varphi(d)
\sum_{\substack{p\le x^{1/3}\\
p\,\equiv\,1\kern-5pt\pmod{du_1u_3}\\(p-1,u_2s)=1}}\frac1{p-1}
\sum_{\substack{q\le x^{1/3}\\
q\,\equiv\,1\kern-5pt\pmod{du_2u_3}\\(q-1,u_1s)=1}}\frac1{q-1}.
\end{equation}
For the sum over $p$ we temporarily ignore the condition that $s\nmid p-1$.  Then $p$
runs over $\varphi_2(u_2)$ residue classes mod $du$.
 In each of these classes, the sum
of $1/(p-1)$ is $\gg 1/\varphi(du)$ by Corollary \ref{cor:apr}.  So, the sum over $p$ appears to
be $\gg \varphi_2(u_2)/\varphi(du)$.  But, we also need to take into account the condition $s\nmid p-1$.
For this, we compute an upper bound for the sum where $s\mid p-1$.
An upper bound sieve result shows that the contribution is $\ll \varphi_2(u_2)/(\varphi(du)\log\log x)$,
which justifies ignoring the condition $s\nmid p-1$.  For the sum over $q$, the analogous
argument shows that it is $\gg \varphi_2(u_1)/\varphi(du)$.

Thus, the expression in \eqref{eq:u1u2u3} is at least of magnitude
\[
\sum_{\substack{d<x^{1/12}\\d\,{\rm squarefree}\\(d,s)=1}}
\frac{u_3\varphi(d)\varphi_2(u_1u_2)}{\varphi(du)^2}\ge
 \sum_{\substack{d<x^{1/12}\\d\,{\rm squarefree}\\(d,s)=1}}\frac{u_3\varphi_2(u_1u_2)}{d\varphi(u)^2}
\gg \frac{u_3\varphi_2(u_1u_2)}{\varphi(u)^2}\log x.
\]
Thus, the proposition now follows from \eqref{eq:m2/3}.
\end{proof}

We are now ready to complete the proof of the lower bound in Theorem \ref{th:3M},
that is, $M_3(x)\gg(\log x)^4$.  From \eqref{eq:m3} and
Proposition \ref{prop:m2u} we have
\[
M_3(x)\gg\sum_{\substack{u\le x^{\delta/12}\\u\,{\rm squarefree}\\s\,\nmid\,u}}\varphi(u)
\sum_{\substack{r\le x^{1/3}\\u\,|\,r-1}}\frac1{r-1}\sum_{\substack{u=u_1u_2u_3\\u_1u_2\,{\rm odd}}}\frac{u_3\varphi_2(u_1u_2)}{\varphi(u)^2}\log x.
\]
It thus follows from Corollary \ref{cor:apr} that
\[
M_3(x)\gg \sum_{\substack{u\le x^{\delta/12}\\u\,{\rm squarefree}\\s\,\nmid\,u}}\sum_{\substack{u=u_1u_2u_3\\u_1u_2\,{\rm odd}}}
\frac{u_3\varphi_2(u_1u_2)}{\varphi(u)^2}\log x.
\]
  We factor the $u$-expression as
\[
\frac{u_3\varphi_2(u_1u_2)}{\varphi(u)^2}=\frac{\varphi_2(u_1)}{\varphi(u_1)^2}
\frac{\varphi_2(u_2)}{\varphi(u_2)^2}\frac{u_3}{\varphi(u_3)^2}.
\]
Note that for $n$ odd we have $\varphi_2(n)/\varphi(n)^2\gg1/n$,
so that
\[
M_3(x)\gg
\log x\Bigg(\sum_{\substack{u_1\le x^{\delta/36}\\u_1\,{\rm odd,\,squarefree}\\s\,\nmid\,u_1}}
\frac{1}{u_1}\Bigg)^2
\sum_{\substack{u_3\le x^{\delta/36}\\u_3\,{\rm squarefree}\\s\,\nmid\,u_3}}\frac{1}{u_3}
\gg(\log x)^4.
\]
This completes the proof of \eqref{eq:3low}.

\section{The lower bound for the second moment}
\label{S:ding}

As mentioned, the proof in Ding \cite{D1} that $M_2(x)\gg\log x$ is not complete
since it relies on an incorrect statement from \cite{MM}.  The proof is easily correctable,
and we give the few details here.  

In view of \eqref{eq:M_k}, we will be done if we show that
\begin{equation}
\label{eq:ding}
\sum_{[p-1,q-1]\le x}\frac1{[p-1,q-1]}\gg\log x.
\end{equation}
Note that the inequality $\sum_{p,q\le x}1/[p-1,q-1]\gg\log x$ is correctly proved in \cite{D1},
and applying this at $\sqrt x$ gives \eqref{eq:ding}.  
We give an alternate proof here.

We may assume that $p\ne q$.
As before,
\[
\frac1{[p-1,q-1]}=\sum_{d\,|\,(p-1,q-1)}\frac{\varphi(d)}{(p-1)(q-1)},
\]
so that
\[
\sum_{[p-1,q-1]\le x}\frac1{[p-1,q-1]}=\sum_{d\le x}\varphi(d)\sum_{\substack{[p-1,q-1]\le x\\d\,|\,(p-1,q-1)}}\frac1{(p-1)(q-1)}.
\]
By placing additional restrictions on $d,p,q$ the expression here only gets smaller.
We do this as follows.  Consider Corollary \ref{cor:apr} from the previous
section with $\varepsilon=1/4$.  We assume that $d$ is squarefree, $d\le x^{\delta/4}$,
 and that $s(x^{1/4})\nmid d$.  We further assume that $p,q\in(x^{1/4},x^{1/2}]$.  
 So, Corollary \ref{cor:apr} implies that $\sum_{p}1/(p-1)\gg1/\varphi(d)$, and the
 same for the sum over $q$.  Thus,
\begin{align*}
\sum_{[p-1,q-1]\le x}\frac1{[p-1,q-1]}
\gg\sum_d\frac1{\varphi(d)}\ge\sum_d\frac1d\gg\log x.
\end{align*}
This completes the proof of \eqref{eq:ding}.

A similar proof can show that $S_2(x)\gg 1$.  Note that the claim that $S_2(x)\asymp 1$
was asserted without proof in \cite{LP}.  Concerning $S_3(x)$, we have a proof that
it is $\gg(\log x)^3$ (and so $S_3(x)\asymp(\log x)^3$ after the result in Section \ref{S:3M}),
but we do not present the details here.

\section{A tail estimate}
\label{S:tail}
In this section we prove the following theorem.
\begin{theorem}
\label{th:lowerbound}
We have
$$
\sum_{\substack{p,\,q\le x\\ [p-1,q-1]> x}}\frac1{[p-1,q-1]}\gg\log x.
$$
\end{theorem}
Recall that in \cite{MM} it is claimed that
\begin{equation}
\label{eq:mm}
\sum_{p,q\le x}\frac1{[p-1,q-1]}=\sum_{[p-1,q-1]\le x}\frac1{[p-1,q-1]}+O(1),
\end{equation}
see the discussion in \cite{MM} at the start of Section 4.  However, the difference between
the two sums in \eqref{eq:mm} is the sum in Theorem \ref{th:lowerbound}, so it cannot
be $O(1)$.

\begin{proof}
We use the full strength of \cite[Theorem 2.1]{AGP} instead of the simplified version used
in Section \ref{sec:low}.
Let $\mathcal D=\mathcal D(x)=\{d\le x^{1/20}:d\hbox{ even},\,\gcd(15,d)=1\}$ with $x$ being sufficiently large. Let $\varepsilon=\delta=.01$, and let $\mathcal D_{\varepsilon,\delta}
=\mathcal D_{\varepsilon,\delta}(x)$
be the possible set of exceptional moduli as described in \cite[Theorem 2.1]{AGP}.
The set $\mathcal D_{\varepsilon,\delta}$ has cardinality $O_{\varepsilon,\delta}(1)$, and
the members are all $>\log x$.  Let
$\mathcal D'=\mathcal D'(x)$ denote the subset of $\mathcal D$ of elements $d$ with $30d$
not divisible by any member of $\mathcal D_{\varepsilon,\delta}(x)$.
 
For each $d\in\mathcal D'$ let $\mathcal P= \mathcal P(x,d)$ denote the set of
primes $p$ with
\begin{itemize}
\item $p\equiv 1\Mod d$,
\item $p\le x$,
\item $\gcd(30,(p-1)/d)=1$.
\end{itemize}
Since $\varphi(30d)/\varphi(d)=16$, it follows from the conditions above that $\mathcal P$ consists of primes $p \le x$ in precisely 3 of the $16\varphi(d)$ reduced residue classes modulo $30d$.  Indeed, if $2^a\,\|\,d$, with $a\ge1$, then 
$p\equiv 2^a+1\Mod{2^{a+1}}$.  Also, $p\equiv2\Mod3$ and $p\equiv2,3,\,\hbox{or}\,4\Mod 5$.

Note that via \cite[Theorem 2.1]{AGP} and partial summation, if $x^{9/10}<t\le x$ with
$x$ sufficiently large, then
\begin{equation}
\label{eq:Pcount}
\Bigg|\sum_{\substack{p\le t\\ p\in\mathcal P}}1-\frac{3t}{16\varphi(d)\log t}\Bigg|\le\frac{6\varepsilon  t}{16\varphi(d)\log t}.
\end{equation}

With $r$ running over primes, let 
\[
f(n)=f(n,x)=\sum_{\substack{7\le r\le x^{1/20}\\r\,|\,n}}\frac1{r-1}.
\]
Note that
\begin{align*}
\sum_{\substack{p\in\mathcal P\\p\le t}}f((p-1)/d)
&=\sum_{7\le r\le x^{1/20}}\sum_{\substack{p\in\mathcal P,\,p\le t\\r\,|\,(p-1)/d}}\frac1{r-1}\\
&<\frac{2t}{(16/3)\varphi(d)\log(t^{8/9}/30)}\sum_{r\ge7}\frac1{(r-1)^2},
\end{align*}
using the explicit version of the Brun--Titchmarsh inequality due to
Montgomery--Vaughan \cite[Theorem 2]{MV}.  Since the final sum here is a constant smaller than $.063$, it follows from \eqref{eq:Pcount} that
\begin{equation}
	\label{eq:sumf}
\sum_{\substack{p\in\mathcal P\\p\le t}}f((p-1)/d)\le \frac3{20}\sum_{\substack{p\in\mathcal P\\p\le t}}1
\end{equation}
for  $x^{9/10}<t\le x$ and $x$ sufficiently large.  Let 
\[
\mathcal P'=\mathcal P'(x,d)=\{p\in\mathcal P:f((p-1)/d)\le1/5\},
\]
so that from \eqref{eq:sumf} we see that
\[
\sum_{\substack{p\in\mathcal P'\\p\le t}}1\ge\frac14\sum_{\substack{p\in\mathcal P\\p\le t}}1
\]
 for $x$ sufficiently large and $x^{9/10}<t\le x$.
Combining this with (\ref{eq:Pcount}) and applying partial summation we obtain
\begin{equation}
\label{eq:P'est}
\sum_{\substack{p\in \mathcal P'\\x^{9/10}<p\le x}}\frac1p > \frac{.0048}{\varphi(d)}\ge\frac{.0096}d
\end{equation}
for all $d\in\mathcal D'$ and $x$ beyond some uniform bound.

For each $d\in\mathcal D'(x)$ let
\[
\mathcal Q=\mathcal Q(x,d)=\{q\le x:q\equiv1\Mod d\},
\]
so that for $x^{9/10}<t\le x$, we have 
\begin{equation}
\label{eq:Qest}
\Bigg|\sum_{\substack{q\in\mathcal Q\\q\le t}}1-\frac t{\varphi(d)\log t}\Bigg|=\left|\pi(t;d,1)-\frac t{\varphi(d)\log t}\right|\le\frac{2\varepsilon t}{\varphi(d)\log t}
\end{equation}
for $x$ sufficiently large.

Next, for $d\in\mathcal D'(x)$ and $p\in\mathcal P'(x,d)$, let 
\[
\mathcal Q'=\mathcal Q'(x,d,p)=\{q\in\mathcal Q: \gcd(q-1,p-1)=d\}.
\]
If $\gcd(q-1,p-1)>d$, then $rd\mid q-1$ for some prime $r\mid (p-1)/d$ with $r\ge 7$
(since $(p-1)/d$ is coprime to 30).  For $x^{9/10}<t\le x$ we have
(using $d\le x^{1/20}$ and $\pi(t;rd,1)\le t/rd$),
\begin{align*}
\sum_{r\,|\,(p-1)/d}\pi(t;rd,1)&=\sum_{\substack{r\,|\,(p-1)/d\\r\le x^{1/20}}}\pi(t;rd,1)
+\sum_{\substack{r\,|\,p-1\\r>x^{1/20}}}\pi(t;rd,1)\\
& \le\sum_{\substack{r\,|\,(p-1)/d\\r\le x^{1/20}}}\frac{2t}{\varphi(d)(r-1)\log(t/rd)}
+\sum_{\substack{r\,|\,p-1\\r>x^{1/20}}}\frac t{rd}\\
& \le\frac94f((p-1)/d)\frac t{\varphi(d)\log t}+O\Big(\frac t{dx^{1/20}}\Big).
\end{align*}
Since $f((p-1)/d)\le 1/5$, we conclude that
\[
\sum_{\substack{q\le t\\q\in\mathcal Q\setminus\mathcal Q'}}1\le\frac {.46t}{\varphi(d)\log t}
\]
for $x$ sufficiently large and $x^{9/10}<t\le x$.  Thus, from \eqref{eq:Qest},
\[
\sum_{\substack{q\le t\\q\in\mathcal Q'}}1\ge\frac t{2\varphi(d)\log t},
\]
so that
\begin{equation}
\label{eq:Q'est}
\sum_{\substack{q\in\mathcal Q'\\ x^{9/10}<q\le x}}\frac1q>\frac{.0525}{\varphi(d)}\ge\frac{.105}d.
\end{equation}

Now for each pair $p,q$ with $p\in\mathcal P'(x,d)$ and $q\in\mathcal Q'(x,d,p)$ with
$p,q>x^{9/10}$,
we have $[p-1,q-1]=(p-1)(q-1)/d>x^{1.75}$.  Further, from \eqref{eq:P'est} and \eqref{eq:Q'est},
\[
\sum_{p,q}\frac1{[p-1,q-1]}=d\sum_p\frac1{p-1}\sum_q\frac1{q-1}>\frac{.001d}{d^2}
=\frac{.001}d.
\]
It remains to note that $\sum_{d\in \mathcal D'}1/d\gg\log x$. In fact, since every member of $ \mathcal D_{\varepsilon,\delta}(x)$ exceeds $\log x$, we have
\[\sum_{a\in \mathcal D_{\varepsilon,\delta}(x)}\sum_{\substack{d\in\mathcal D\\a\mid 30d}}\frac{1}{d}=O_{\varepsilon,\delta}(1),\]
so that
\[\sum_{\substack{p,\,q\le x\\ [p-1,q-1]> x}}\frac1{[p-1,q-1]}>.001\sum_{d\in \mathcal D}\frac{1}{d}+O(1)=\frac{1}{75000}\log x+O(1)\]
for all sufficiently large $x$. This completes the proof.
\end{proof}

Let
\[
\mathcal M_2(x)=\sum_{[p-1,q-1]\le x}\frac1{[p-1,q-1]},\quad
\mathcal M'_2(x)=\sum_{p,q\le x}\frac1{[p-1,q-1]}.
\]
Note that $\mathcal M_2(x)=M_2(x)+O(1)$.\footnote{~In fact, we have
$0\le\mathcal M_2(x)  -M_2(x)\le S_2(x)$.  If we had $S_2(x)=C+o(1)$, then we might
expect that $(\mathcal M_2(x)-M_2(x))/S_2(x)=1-\gamma+o(1)$, where $\gamma$ is
the Euler--Mascheroni constant.  The numbers in Table \ref{Ta:xM_2(x)&xS_2(x)} and Table \ref{Ta:MM&D}
strongly support such a relation.}
The sum $\mathcal M'_2(x)$ is heuristically shown to be
$\sim 2\zeta(2)\zeta(3)/\zeta(6)\log x$ in \cite{D2}.  We have shown above that $\mathcal M'_2(x)-\mathcal M_2(x)\gg\log x$.
It would be interesting to either prove or give a heuristic for $\mathcal M'_2(x)-\mathcal M_2(x)\sim\kappa\log x$
for some explicit number $\kappa>0$.  In lieu of this we provide some numerical experiments.
The difference $\mathcal M'_2(x)-\mathcal M_2(x)$ does indeed look like it is growing linearly
in $\log x$ with a slope of about $0.69$.

\begin{table}[ht]
		\caption{Numerical values of $\mathcal M_2(10^k)$ and $\mathcal M'_2(10^k)$}
		\label{Ta:MM&D}
		\begin{tabular}{|l|l|l|c|c|} \hline
			~$k$  & $\mathcal M_2(10^k)$ &$\mathcal M'_2(10^k)$&$\mathcal M'_2(10^k)-\mathcal M_2(10^k)$&$0.69\log(10^k)-2.7$\\ \hline
			~$3$ & 16.6272 &19.0012&2.3740 &2.0664\\
			~$4$ & 23.0838&26.9182& 3.8347&3.6551\\
			~$5$ & 29.7107&35.0582& 5.3475&5.2439\\
			~$6$ & 36.5061&43.3902&6.8841&6.8327\\
			~$7$ & 43.3932&51.8341& 8.4409&8.4214\\
			~$8$ & 50.3485&60.3521&10.0036 &10.0103\\
			~$9$ & 57.3533&68.9220&11.5687&11.5991\\
						\hline
		\end{tabular}
	\end{table}

\section*{Acknowledgments}
We thank Nathan McNew, Paul Pollack, and the referee for some helpful comments.

\end{document}